\documentclass{article}

\begin{document}
\begin{center}
{\bf\large On Shanks' Algorithm for Modular Square Roots}
\end{center}
\begin{abstract}
Let $p$ be a prime number, $p=2^nq+1$, where $q$ is odd. D. Shanks
described an algorithm to compute square roots $\pmod{p}$ which needs
$O(\log q + n^2)$ modular multiplications. In this note we describe two
modifications of this algorithm. The first needs only $O(\log q + n^{3/2})$
modular multiplications, while the second is a parallel algorithm
which needs $n$ processors and takes $O(\log q+n)$ time.
\end{abstract}

MSC-Index 11Y16, 68Q25, 68W10


Short title: Computing Modular Square Roots

In \cite{Sha}, D. Shanks gave an efficient algorithm for computing
square roots modulo a prime. If $p=2^nq+1$, this algorithm consists of
an initialization, which takes $O(\log q)$ modular multiplications,
and a loop, which is performed at most $n$ times and needs $n$ modular
multiplications at most. Hence the total cost are $O(\log q + n^2)$
modular multiplications. This is actually the normal running time, for
S. Lindhurst \cite{Lin} has shown that on average the loop needs
$\frac{1}{4}(n^2+7n-12) + 1/2^{n-1}$ modular multiplications. For most
prime numbers $p$, $n$ is much smaller then $\sqrt{\log q}$, hence the
initialization will be the most costly part, 
however, prime numbers occuring ``in practice'' are not necessarily
random, and if $p-1$ is divisible by a large power of 2, the loop
becomes more expensive then the initialization. In this note we
will give two modifications of Shanks' algorithm. The first
algorithm needs only $O(\log q + n^{3/2})$ modular
multiplications, while the second is a parallel algorithm running on
$n$ processors which needs $O(\log q + n)$ time. On the other
hand both our algorithms have larger space 
requirements. Whereas Shanks' algorithm has to store only a bounded
number of residues $\pmod{p}$, our algorithms have to create two fields,
each containing $n$ residues $\pmod{p}$. However, on current hardware
this amount of memory appears easily manageable compared to the 
expenses of the computation.

We assume that looking up an element in a table of
length $n$ is at most as expensive as a modular
multiplication, an assumption which is certainly satisfied on any
reasonable computer.

First we give a description of Shanks' algorithm.
We assume that we are given a prime $p=2^nq+1$, a quadratic residue
$a$ and a noresidue $n$, and are to compute an $x$ such that
$x^2\equiv a\pmod{p}$. Then the algorithm runs as follows.
{\bf Algorithm 1:\\}
\begin{enumerate}
\item Set $k=n$, $z=u^q$, $x=a^{(q+1)/2}$, $b=a^q$.
\item Let $m$ be the least integer with $b^{2^m}\equiv 1\pmod{p}$.
\item Set $t=z^{2^{k-m-1}}$, $z=t^2$, $b=bz$, $x=xt$.
\item If $b=1$, stop and return $x$, otherwise set $k=m$ and go to
step 2.
\end{enumerate}

It is easy to see that the congruence $x^2\equiv ab\pmod{p}$ holds at
every stage of the algorithm, hence, if it terminates we really
obtain a square root of $a$.

To show that this algorithm terminates after at most $n$ loops,
consider the order of $b$ and $z\pmod{p}$. After the first step, the
latter is $2^n=2^k$, since $u$ is a nonresidue, whereas the first one is
strictly smaller, since $a$ is a quadratic residue. In the second step
the order of $b$ is determined to be exactly $2^m$, and in the third
step $z$ is replaced by some power, such that the new value of $z$ has
order exactly $2^m$, too. Then $b$ is replaced by $bz$, thus the order
of the new value of $b$ is $2^{m-1}$ at most. Setting $k=m$, we get
the same situation as before: the order of $z$ is exactly $2^k$, and
the order of $b$ is less. Hence every time the loop is executed, 
the order of $b$ is reduced, at the same time it always remains a
power of 2. Hence after at most 
$n$ loops, the order of $b$ has to be 1, i.e. $b\equiv 1\pmod{p}$.

The next algorithm is our first modification of Algorithm 1.\\[2mm]
{\bf Algorithm 2:}
\begin{enumerate}
\item Set $k=n$, $z=u^q$, $x=a^{(q+1)/2}$, $b=a^q$. 
\item Compute $z^2, z^{2^2}, z^{2^3}, \ldots, z^{2^n}$ and store these
values in an array.
\item Compute $b^2, b^{2^2}, b^{2^3}, \ldots, b^{2^n}$ and store these
values in an array.
\item Set $i=1$, $b_0=b$, $z_0=z$
\item Let $m$ be the least integer, such that $b_0^{2^m}z_1^{2^m}\cdots
z_i^{2^m}\equiv 1\pmod{p}$.
\item Set $t=z_i^{2^{k-m-1}}$, $z_{i+1}=t^2$, $b=bz_{i+1}$, $x=xt$,
$i=i+1$, $k=m$.
\item If $b=1$, stop and return $x$.
\item If $i<\sqrt{n}$, continue with 5, otherwise set $z=z_{i+1}$ and
continue with 3.
\end{enumerate}

First observe that there are no essential changes to the
algorithm. The only 
difference is that in step 5 - which corresponds to step 2 in the original
algorithm - no explicite reference to $b$ is made, but $b$ is replaced
by $b_0z_1\cdots z_i$. Of course, the numerical value of these
expressions is the same, however, we claim
that in the form above the algorithm needs only $O(\log q+n^{3/2})$
modular multiplications.

Note first that for any $i$ at any stage in the algorithm,
$z_i=u^{q\cdot 2^l}$ for some integer $l$, and the same is true for
$t$. In fact, the only point where some operations are performed with
these numbers is in line 6, where a certain number of squarings are
performed, however, the effect of this operation is just a shift
within the array of precomputed values. Hence, for any
exponent $m$ and index $i$, $z_i^{2^m}$ can be obtained by looking up
in the array generated in step 2. After this remark we can compute the
running time. The inner loop is performed at most $n$ times, hence
step 6 needs $O(n)$ modular multiplications alltogether. The outer
loop is performed at most $[\sqrt{n}]$-times, hence step 3 requires
$n^{3/2}$ modular multiplications alltogether. Step 2 requires n
multiplications and is performed once, and steps 1, 4, 7 and 8 can
be neglected.

Hence we have to consider step 5. The check whether for a
given $m'$ the congruence $b^{2^{m'}}z_1^{2^{m'}}\cdots z_i^{2^{m'}}\equiv
1\pmod{p}$ holds true, can be done using $i$ modular multiplications,
since all the 
powers can be obtained by looking up in the arrays generated in step 2
and 3. We already know at this stage that the congruence holds for
$m'=k$, hence we compute the product for $m'=k-1,\, k-2,\, \ldots$, untill
we find a value for $m$ such that the product is not
$1\pmod{p}$. Doing so we have to check $k-m$ values $m'$, hence at a 
given stage this needs $(k-m)i=O((k-m)\sqrt{n})$ modular
multiplications. To estimate the sum of these costs, introduce a
counter $\nu$, which is initialized to be $0$ in step 1 and raised by
one in step 5, that is, $\nu$ counts the number of times the inner
loop is executed. Define a sequence $(m_n)$, where $m_n$
be the value of $m$ as found in step 5 when $n=\nu$. With this
notation the costs of step 5 as estimated above are
$O((m_{\nu-1}-m_\nu)\sqrt{n})$, and the sum over $\nu$ telescopes. Since
$m_1\leq n$, and $m_{\nu_1}=1$, where $\nu_1$ is the value of the
counter $\nu$ when the algorithm terminates, the total cost of step 5
is $O(n^{3/2})$.

Putting the estimates together we see that there is a total amount of
$O(n^{3/2})$ modular multiplications. In the same way one sees that we
need $O(n^{3/2})$ look ups, and by our assumption on the costs of the
latter operation we conclude that the running time of Algorithm 2 is
indeed $O(\log q + n^{3/2})$.

Finally we describe a parallel version of Algorithm 1:\\[2mm]
{\bf Algorithm 3:}
\begin{enumerate}
\item Set $k=n$, $z=u^q$, $x=a^{(q+1)/2}$, $b=a^q$. 
\item Compute $z^2, z^{2^2}, z^{2^3}, \ldots, z^{2^n}$ and store these
values in an array.
\item Compute $b^2, b^{2^2}, b^{2^3}, \ldots, b^{2^n}$ and store these
values in an array.
\item Let $m$ be the least integer, such that $b^{2^m}\equiv 1\pmod{p}$.
\item Set $t=z^{2^{k-m-1}}$, $z=t^2$, $x=xt$, $k=m$.
\item Set $b=bz$, compute $b^2, b^{2^2}, \ldots, b^{2^m}$ and replace
the powers of $b$ by these new values.
\item If $b=1$, stop and return $x$, otherwise continue with step 4.
\end{enumerate}

It is clear that this algorithm is equivalent to Algorithm 1,
furthermore all steps with the exception of step 6 can be performed by
a single processor in time $O(\log q + n)$. Now consider step 6. This
step has to be executed at most $n$ times, and we claim that it can be
done by $n$ processors in a single step. Indeed, since all relevant
powers both of the old value of $b$ and of $z$ are stored, each of 
the powers of the new value of $b$ can be obtained by a single
multiplication, and all these multiplications can be done
independently from each other on different processors. Hence,
Algorithm 3 runs in time $O(\log q + n)$ on $m$ processors.

Jan-Christoph Puchta\\
Mathematical Institute\\
24-29 St Giles\\
Oxford OX1 3LB\\
United Kingdom\\
puchta@maths.ox.ac.uk\\[4mm]

\begin{thebibliography}{5}
\bibitem{Lin} J. Lindhurst, {\em An analysis of Shanks's algorithm for
computing square roots in finite fields}, in: Number theory
(Ottawa, 1996), CRM Proc. Lecture Notes, 19, 1999, 231--242 
\bibitem{Sha} D. Shanks, {\em Five number-theoretic algorithms}, in:
Proceedings of the Second Manitoba Conference on Numerical
Mathematics, Congressus Numerantium, No. VII, 1973, 51--70 
\end{thebibliography}
\end{document}